\documentclass[a4paper,11pt,oneside]{amsart}
\usepackage{amsbsy,amssymb,pstricks,pst-node,mathrsfs,MnSymbol}
\usepackage{mathtools,graphicx}
\usepackage{subfigure}
\usepackage[shortlabels]{enumitem}
\usepackage[pagewise]{lineno}
\usepackage{tikz}
\usetikzlibrary{positioning}
\newtheorem{theorem}{Theorem}[section]

\newtheorem{lemma}[theorem]{Lemma}
\newtheorem{remark}[theorem]{Remark}

\theoremstyle{case}

\newtheorem{corollary}[theorem]{Corollary}
\newtheorem{example}[theorem]{Example}
\newtheorem{note}[theorem]{Note}
\theoremstyle{definition}
\newtheorem{definition}[theorem]{Definition}

\newcommand{\covering}{%
	\mathrel{-\mkern-4mu}<%
}

\begin{document}
	

	%
	\title{On the strong metric dimension of the zero-divisor graph of a lattice  }
\maketitle
\markboth{Pravin Gadge and Vinayak Joshi}{On the strong metric dimension of  the zero-divisor graph of a lattice  }\begin{center}\begin{large}Pravin Gadge$^\text{a}$ and  Vinayak Joshi$^\text{b}$\end{large}\\\begin{small}\vskip.1in$^\text{a}$\emph{GES's Shri Bhausaheb Vartak Arts, Commerce and Science College, Borivali - 400091, Maharashtra, India.}\\$^\text{b}$\emph{Department of Mathematics, Savitribai Phule Pune University, Pune - 411007, Maharashtra, India.\\E-mail: praving2390@gmail.com (P. Gadge),  vvjoshi@unipune.ac.in (V. Joshi) }\end{small}\end{center}\vskip.2in

\begin{abstract} In this paper, the generalized blow-up of a Boolean lattice  $L\cong \textbf{2}^n$ using finite chains is introduced. Also, we compute the strong metric dimension of the zero-divisor graph of the blow-up of a Boolean lattice. These results are applied to calculate the strong metric dimension of the comaximal graph, the comaximal ideal graph, the zero-divisor graph of a reduced ring, and the component graph of a vector space.  	
\end{abstract}\vskip.2in
\noindent\begin{Small}\textbf{Mathematics Subject Classification (2020)}:
	05C25, 06A07, 05C17, 13A70.    \end{Small}\vskip.2in
\vskip.25in
\baselineskip 17truept 

Beck \cite{Be} originally introduced the idea of associating a graph with a commutative ring, primarily focusing on colorings. Anderson and Livingston \cite{AL} subsequently modified the definition of the zero-divisor graph of a commutative ring $R$, denoted by $\Gamma(R)$. In this graph, the vertex set is the set of all nonzero zero-divisors of $R$, and two vertices $x$ and $y$ are adjacent if $xy=0$. Many researchers have studied the interplay between the ring-theoretic properties of $R$ and the graph-theoretic properties of $\Gamma(R)$.

Harary and Melter \cite{hara} first introduced the concept of the metric dimension of a graph. In 2004, Seb\"{o} and Tannier \cite{sebo} introduced a more specific parameter known as the strong metric dimension. Many researchers have studied the concept of metric dimension and the strong metric dimension for a wide variety of graphs, such as Cayley graphs, trees and unicyclic graphs, wheel graphs, Cartesian product graphs, etc. (see \cite{gcl}, \cite{skb}, \cite{galai}).

Finding the metric and the strong metric dimensions of graphs are NP-complete problems. Thus, some researchers have been interested in determining these parameters for graphs of algebraic structures and ordered structures; see \cite{ip}, \cite{dk}, and \cite{galai}.

In this paper, the generalized blow-up of a Boolean lattice $L \cong \textbf{2}^n$ using finite chains is introduced. Additionally, we compute the strong metric dimension of the zero-divisor graph of the blow-up of a Boolean lattice. These results are applied to calculate the strong metric dimension of the comaximal graph, the comaximal ideal graph, the zero-divisor graph of a reduced ring, and the component graph of a vector space.

\section{Preliminaries}

By $G = G(V, E)$, we mean a simple and undirected graph $G$ with the vertex set $V = V(G)$ and the edge set $E = E(G)$. Let $N(v)$ denotes the set of all vertices adjacent to a vertex $v$ in $G$, and $N[v] = N(v) \cup \{v\}$. A set $S$ of vertices in $G$ forms a \textit{vertex cover} if every edge of $G$ has at least one end in $S$. The \textit{vertex cover number} of $G$, denoted by $\alpha(G)$, is the minimum cardinality required for a vertex cover of $G$. An \textit{independent set} of a graph $G$ is a set of vertices such that no two vertices are adjacent. The \textit{independence number} of $G$, denoted by $\beta(G)$, is the cardinality of a largest independent set in $G$.

For a connected graph $G$, consider a subset $S = \{v_1, v_2, \dots, v_k\}$ of $V(G)$, and let\linebreak  $v \in V(G) \setminus S$. The \textit{metric representation} of $v$ with respect to $S$ is expressed as the \linebreak $k$-vector (ordered $k$-tuple) $D(v|S) = (d(v, v_1), d(v, v_2), \dots, d(v, v_k))$. If, for $S \subseteq V(G)$, the equality $D(u|S) = D(v|S)$ holds for every pair of $u, v \in V(G) \setminus S$, implying $u = v$, then $S$ is referred to as a \textit{resolving set} for $G$. The metric basis for $G$ is a resolving set $S$ with the minimum cardinality, and the number of elements in $S$ is defined as the \textit{metric dimension} of $G$, denoted by $dim_M(G)$.

In a connected graph $G$, a vertex $w$ is said to \textit{strongly resolve} two vertices $u, v$, if there exists a shortest path from $u$ to $w$ containing $v$ or a shortest path from $v$ to $w$ containing $u$. A set $W$ of vertices is termed a \textit{strong resolving set} for $G$, if every pair of vertices in $G$ is strongly resolved by at least one vertex in $W$. The smallest cardinality of a strong resolving set for $G$ is named the \textit{strong metric dimension} of $G$, denoted by, $sdim_M(G)$.

\indent Let $P$ be a \textit{partially ordered set} (poset) with $0$. 
 Given any $A \subseteq P$, the \textit{upper cone} of $A$ is the set $A^u=\{b \in P \mid a \leq b$ for every $a \in A\}$ and the \textit{lower cone} of $A$ is the set $A^{\ell}=\{b \in P \mid b \leq a$ for every $a \in A\}$.
 The \textit{annihilator} of $A$, denoted by $A^\perp$, is the set of elements $b$ in $P$ such that $\{a, b\}^{\ell} = \{0\}$ for all $a \in A$. If $A = \{a\}$, then $A^\perp$ is denoted by $a^\perp$. Let $(P,\leq)$ be a poset, then  the \textit{dual} of $P$ is denoted by $(P^{\partial}, \geq)$ is the poset with the partial order $a \geq b$ in $P^{\partial}$ if and only if $a \leq b$ in $P$.
 
  Let $x$ and $y$ be elements of $P$. Then $y$ \textit{covers} $x$, written $x\covering y$, if $x < y$ and there is no element $z$ such that $x < z < y$. If $0 \covering x$, then $x$ is called an \textit{atom} of $P$. Moreover, $P$ is called \textit{atomic} if every nonzero element contains an atom. The set of atoms of $P$ is denoted by $Atoms(P)$. By a \textit{chain}, we mean a poset in which any two elements are comparable. If $a$ and $b$ are incomparable elements of $P$, then we denote it by $a || b$.

\indent A poset $P$ is said to be \textit{bounded}, if $P$ has both the least element $0$ and the greatest element $1$. An element $b$ of a bounded poset $P$ is a \textit{complement} of $a \in P$ if $\{a, b\}^{\ell} = \{0\}$ and $\{a, b\}^{u} = \{1\}$. A \textit{pseudocomplement} of $a \in P$ is an element $b \in P$ such that $\{a, b\}^{\ell} = \{0\}$, and if $\{a, x\}^{\ell} = \{0\}$, then $x \leq b$. It is easy to confirm that for any element $a$ in $P$, there is at most one pseudocomplement, denoted as $a^{\ast}$ if it exists. A bounded poset $P$ is called \textit{complemented} (respectively, \textit{pseudocomplemented}) if every element of $P$ has a complement (respectively, $a^{\ast}$ exists for every $a \in P$). A bounded poset $P \cong M_n$ if and only if every element of $P \setminus \{0,1\}$ is an atom, if and only if every element of $P \setminus \{0,1\}$ is covered by $1$.

Define a \textit{zero-divisor} of $P$ to be any element of the set $Z(P) = \{a \in P \mid \text{there exists } b \in P \setminus \{0\} \text{ such that } \{a, b\}^{\ell} = \{0\}\}$. An element $a \in P$ is called \textit{dense} if $a \notin Z(P)$. The set of all dense elements of $P$ is denoted by $D(P)$. As in \cite{tuwu}, the zero-divisor graph of $P$ is the graph $G(P)$ whose vertices are the elements of $Z^*(P) = Z(P) \setminus \{0\}$ such that two vertices $a$ and $b$ are adjacent if and only if $\{a, b\}^{\ell} = \{0\}$.

\par   Let $a$ be any element of a lattice $L$. The ideal generated by $a$ is called {\it principal ideal}. It is denoted by $(a]$ and is given by $(a]=\{x\in L \mid x\leq a\}$. Dually, we have the concept of a {\it principal filter}. A lattice $L$ is called a \textit{$0$-distributive lattice}, \index{$0$-distributive lattice}  
if $a\wedge b=0$ and $a\wedge c=0$ implies $a\wedge(b\vee c)=0$. Dually, we have the concept of a  \textit{$1$-distributive lattice}. Moreover, a bounded  distributive  and complemented poset $P$
is called \emph{Boolean}. It is well-known that in a Boolean lattice, complementation coincides with pseudocomplementation (cf. \cite[Lemma 2.4]{vjkh}). In particular, if $P$ is Boolean, then $P$ is pseudocomplemented, and every element $x\in P$ has the unique complement $x'$. Sometimes, it is also denoted by $x^\ast$.

\section{\bf Strong Metric Dimension of  the Zero-Divisor Graph of  a Lattice}
Through a series of papers (see \cite{sdj1, pgg, nkvj23, ncv}), it has been observed that the zero-divisor graph of ordered sets serves as a tool to study various graphs associated with algebraic structures. Notably, the blow-up of a Boolean lattice acts as a prototype for studying these graphs, such as the comaximal graph of a ring, the nonzero component graph of a vector space, and the zero-divisor graph of a reduced ring. 

In this section, we derive a formula for the strong metric dimension of the zero-divisor graph of a blow-up of a Boolean lattice. As a result, this formula also applies to the strong metric dimension of the aforementioned graphs.

In the existing literature, researchers have used two graphs, $G_{SR}$ and $G^{**}$, to determine the strong metric dimension of a graph $G$. By Theorem \ref{gala}, finding the strong metric dimension of $G$ requires determining the vertex cover of $G_{SR}$. However, the structure of $G_{SR}$ is somewhat complex. Hence, a new graph, $G^{**}$, is introduced, and its relationship with $G_{SR}$ is explored. Consequently, the problem of finding the strong metric dimension of $G$ reduces to finding the vertex cover number of $G^{**}$.

	We start by introducing the requisite background definitions and findings.

\begin{theorem}[{D. Lu and T. Wu \cite[Proposition 2.1]{tuwu}}] \label{diameter}
Let $P$ be a poset. Then $G(P)$ is connected graph with diam$(G(P))\leq 3$.
\end{theorem}

\begin{lemma}\label{finite}
\par Let $P$ be a poset with $0$. Then $dim_M (G(P))$ is finite if and only if $G(P)$ is finite. 
\end{lemma}  
\begin{proof}
Assume that  $dim_M (G(P))$ is finite. Let $W$ be the metric basis for $G(P)$ with $|W|=k$ for some non-negative integer $k$. By Theorem \ref{diameter}, the diameter of $G(P)$ is at most $3$, i.e., $d(x,y)\in \{1,2,3\}$ for every distinct $x,y\in V(G(P))$. Then for each $x\in V(G(P))$, the metric representation $D(x|W)$ is the $k$-coordinate vector, where each coordinate is in the set $\{1,2,3\}$. Thus, there are only $3^k$ possibilities for $D(x|W)$. Since $D(x|W)$ is unique for each $x\in V(G(P))$, so $|V(G(P))|\leq 3^k$. This implies that $V(G(P))$ is finite. Hence, $G(P)$ is finite.
The converse is obvious.	
\end{proof}

It is easy to observe that every strong resolving set of a graph $G$ is also a resolving set. Hence $\operatorname{dim}_M(G) \leq \operatorname{sdim}_M(G)$.

\begin{corollary}
Let $P$ be a poset. Then $\operatorname{sdim}_M(G(P))$ is finite if and only if $G(P)$ is finite.
\end{corollary}

A widely recognized result, credited to Gallai, establishes the connection between the independence number $\beta(G)$ and the vertex cover number $\alpha(G)$ of a graph $G$.

\begin{theorem}[Gallai's Theorem] \label{galai}
For any graph $G$ of order $n$, $\alpha(G) + \beta(G) = n$.
\end{theorem}

\begin{definition}
A vertex $u$ in a graph $G$ is considered {\it maximally distant} from $v$, if for every $w$ in the neighborhood of $u$, the distance from $v$ to $w$ is less than or equal to the distance from $u$ to $v$. When both $u$ is maximally distant from $v$ and $v$ is maximally distant from $u$, we describe $u$ and $v$ as mutually maximally distant.	\end{definition} 
Note that if $u$ is maximally distant from $v$, then  $v$ need not be maximally distant from $u$. Also, a vertex $u$ is not maximally distant from itself. One can see that in a graph $G(P)$ shown in Figure \ref{figure1}, the vertex $(1,0,1)$ is maximally distant from $(0,1,0)$, however $(0,1,0)$ is not maximally distant from $(1,0,1)$.
\par The {\it boundary} of $G$, denoted by $\partial(G)$, consists of vertices $u$ in $V(G)$ for which there exists a vertex $v$ in $V(G)$ such that $u$ and $v$ are mutually maximally distant.
The boundary of the graph $G(L)$ is given in Example \ref{example1} (2).

The concept of a strong resolving graph was  introduced in  \cite{galai}.
\begin{definition} [Oellermann and Peters-Fransen \cite{galai}]
Let $G$ be a graph. The \textit{strong resolving graph} of $G$, denoted by $G_{SR}$, with the vertex set  $\partial(G)$ and two distinct vertices $u$ and $v$ are connected in $G_{SR}$ if and only if $u$ and $v$ are mutually maximally distant in $G$.
\end{definition}

It was proved in \cite[Theorem 2.1]{galai} that determining the strong metric dimension of a graph $G$ is nothing but the vertex cover number of $G_{SR}$.

\begin{theorem}[{Oellermann and Peters-Fransen \cite[Theorem 2.1]{galai}}] \label{gala}
For any connected graph $G$, $sdim_M(G) = \alpha(G_{SR})$.
\end{theorem}

\begin{example}\label{example1}
\begin{enumerate}
\item  Let $L=M_{n}$, then $G(L)=K_n$. Since $(K_n)_{SR}=K_n$, we have $sdim_M(G(L))=n-1$.
\item  Let $L=C_2\times C_2\times C_2$ and $G(L)$ be its zero-divisor graph. Suppose that  $X = \{(1, 0, 0), (0,1, 0), (0, 0,1)\}$
and $Y =\{(1,1, 0), (0,1,1), (1, 0,1)\} $. One can see that for any $u \in X$, there is no $v \in V(G(L))$ such that $u$ and $v$ are mutually maximally distant, whereas each pair of vertices in $Y$ are mutually maximally distant. This gives that $\partial(G(L))=\{(1,1, 0), (0,1,1), (1, 0,1)\}$ and $G(L)_{SR} = K_{3}$. Since $\alpha(G(L)_{SR})=2$, by Theorem \ref{gala},  $sdim_{M}(G(L))=2$. On the other hand, $W = \{(1, 1, 0), (0,1, 1)\}$ is a
minimum  cardinality strong resolving set, i.e., $sdim_M(G(L)) = 2$. Note that the strong metric dimension of a graph isomorphic to $G(L)$ is calculated in \cite{smco}.

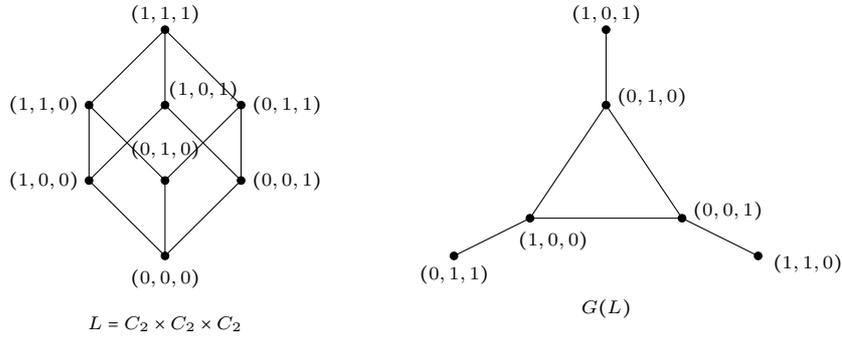
\begin{figure}[h]

	\begin{center}
		\begin{tikzpicture}	[scale=1]		
			\begin{scope}[shift={(-1,0)}]
				\draw [fill=black] (0,0) circle (0.05);
				\draw [fill=black] (-1,1) circle (0.05);
				\draw [fill=black] (0,1) circle (0.05);
				\draw [fill=black] (1,1) circle (0.05);
				\draw [fill=black] (-1,2) circle (0.05);
				\draw [fill=black] (0,2) circle (0.05);
				\draw [fill=black] (1,2) circle (0.05);
				\draw [fill=black] (0,3) circle (0.05);

				\draw (0,0) -- (-1,1) -- (-1,2) -- (0,3);
				\draw (0,0) -- (0,1) -- (-1,2) ;
				\draw (0,0) --  (1,1) -- (1,2) -- (0,3);
				\draw (-1,1) -- (0,2) -- (0,3);
				\draw (0,1) -- (1,2);
				\draw (1,1) -- (0,2);
				\node at (0,-0.9) {\tiny  $ L=C_2\times C_2 \times
					C_2$};
				
				\node at (0,-0.3) {\tiny $(0,0,0)$};
				\node at (-1.6,1) {\tiny $(1,0,0)$};
				\node at (-1.6,2) {\tiny $(1,1,0)$};
				\node at (0,3.2) {\tiny $(1,1,1)$};
				\node at (0,1.4) {\tiny $(0,1,0)$};
				\node at (1.6,1) {\tiny $(0,0,1)$};
				\node at (1.6,2) {\tiny $(0,1,1)$};
				\node at (0.5,2.2) {\tiny $(1,0,1)$};

			\end{scope}

			\begin{scope}[shift={(0.8,0)}]

				\draw [fill=black] (4,3) circle (0.05);
				\draw [fill=black] (5,0.5) circle (0.05);
				\draw [fill=black] (3,0.5) circle (0.05);
				\draw [fill=black] (4,2) circle (0.05);
				\draw [fill=black] (6,0) circle (0.05);
				\draw [fill=black] (2,0) circle (0.05);

				\draw (3,0.5) -- (4,2) -- (4,2)  -- (5,0.5);
				\draw (2,0)--(3,0.5);
				\draw (3,0.5) -- (5,0.5);
				\draw (6,0)--(5,0.5);
				\draw (4,3)--(4,2);
				
				\node at (3.3,0.2) {\tiny  $(1,0,0)$};
				\node at (4,3.2) {\tiny $(1,0,1)$};
				\node at (5.6,0.6) {\tiny $(0,0,1)$};
				\node at (6.65,-0.1
				) {\tiny $(1,1,0)$};
				\node at (2,-0.25) {\tiny $(0,1,1)$};
				\node at (4.6,2.1) {\tiny $(0,1,0)$};
				\node at (4,-0.7) {\tiny $G(L)$};

			\end{scope}

		\end{tikzpicture}
		\caption{A Boolean lattice $L$ and its zero-divisor graph  $ G(L)$}\label{figure1}
	\end{center}
\end{figure}

\item The zero-divisor graph $G(L^B)$ of the lattice $L^B$ is shown in Figure \ref{figure4}. The  set   $W=\{x_{1}^{1},x_{1}^{2},x_{3}^{1},x_{13}^{1},x_{13}^{2},x_{13}^{3},x_{12}^{1},x_{12}^{2}\}$ is a minimum cardinality strong resolving set for $G(L^B)$. Hence $sdim_{M}(G(L^B))=8$.
\end{enumerate}
\end{example}

\begin{definition}[\cite{jwp}]\label{defn1}
Let $L$ be a lattice with $0$. Define a relation $\sim$ on $L$ as $x\sim y$ if and only if $x^\perp=y^\perp$. Clearly, $\sim$ is an equivalence relation on $L$. Let  $[a]$ denotes the equivalence class of $a$ under $\sim$. The set of equivalence classes of $L$ will be denoted by $[L]$=\{$[a]\mid a \in L$\}. Note that $[L]$ is a meet-semilattice under the partial order given by ~$[a]\leq [b]$ if and only if $b^\perp\subseteq a^\perp$ with $[a]\wedge [b]=[a\wedge b]$ $($see {\cite[Lemma 2]{jwp}}$)$ . If $L$ is 0-distributive, then $[L]$ is a lattice; see \cite{jwp}. 
This result is further strengthened by Khandekar and Joshi \cite{nkvj23}. 
\end{definition}

\begin{theorem} [{Khandekar and Joshi \cite[Theorem 1.1]{nkvj23}}]\label{paka}
Let $L$ be a $0$-distributive bounded lattice with finitely many atoms. Then $[L]$ is a Boolean lattice.
\end{theorem}

The blow-up of a graph was first introduced by M. Ye et al. in \cite{mytwl}. On similar lines, Gadge et al. \cite{pgg} defined the blow-up of a Boolean lattice using finite chains with a certain total order. Now, we introduce the generalized blow-up of a Boolean lattice.

\textbf{Throughout the paper, let $L$ be a Boolean lattice with atoms $q_i$, $1\leq i\leq n$ ($n \geq 3$), i.e., $L\cong \mathbf{2}^n $ with $|L|\geq 8$. }

It is known that every element of $L$ is the join of atoms in $L$. Hence an element $x=q_{i_1}\vee q_{i_2}\vee \dots \vee q_{i_k}\in L$ ($\{i_1,i_2,\dots,i_k\}$$\subseteq$ $\{1,2,\dots,n\}$) can also be represented  as $(x_1,x_2,\dots,x_n)\in \mathbf{2}^n$ where $$
x_j=\begin{cases}
1 & \text{if $j\in \{i_1,i_2,\dots , i_k\}$ }\\
0 & \text{otherwise}.
\end{cases}
$$\\
Thus,  an atom $q_i\in L$ is denoted by $(0,\dots,0,1,0,\dots,0)$, where $1$ is at $i^{th}$ position.

\begin{definition}\label{blowup}
The blow-up $L^B$ of a Boolean lattice $L\cong \mathbf{2}^n$ using chains
is   obtained as follows:
\begin{enumerate}
\item Replace each atom $q_i$ ($1\leq i\leq n$) of $L$ by a chain $C_i$ of finite  length, say $m_i-1$, with elements $q_i=x_i^{1},x_i^{2},\dots,x_i^{m_i}$ such that  $x_{i}^1 \covering x_{i}^2 \covering \dots \covering x_{i}^{m_i}  $. 
\item Let $ x=\bigvee\limits_{j=1}^{k}q_{i_j}\in L\setminus\{1\}$, where $q_{i_j}$ be atoms of $L$ with $i_j\in\{1,2,\dots, n\}$. Replace $x\in L$ by a chain 		
$C_{i_1i_2 \dots i_k}$ of  finite length, say  $n_{j}-1$,  with elements $x=x_{i_1i_2 \dots i_k}^{1}, ~ x_{i_1i_2 \dots i_k}^{2},~ \dots, x_{i_1i_2 \dots i_k}^{n_j}$ for some $n_{j}\in \mathbb{N}$  such that 
$x_{i_1i_2 \dots i_k}^{1}\covering x_{i_1i_2 \dots i_k}^{2}\covering \dots \covering x_{i_1i_2 \dots i_k}^{n_{j}}$,
where $\{i_1,i_2,\dots, i_k\} \subseteq \{1,2,\dots, n\} $.
\item The elements $0$ and  $1$ of $L$ will be represented by  $\mathbf{0}$ and $\mathbf{1}$ in $L^B$ respectively.

\end{enumerate}
\end{definition}
\begin{remark}\label{diam}
Note that diam$(G(L^B))=3=\text{diam}(G(L))$, as $L\cong \mathbf{2}^n $ with $n \geq 3$.
\end{remark} 
We will represent the elements of $L^B$ in terms of tuples as follows.\\
An element $x^{t}_{i_1i_2\dots i_k}$  ($1\leq t \leq n_j$ for some $n_j\in \mathbb{N}$) on the chain $C_{i_1i_2\dots i_k}$ ($\{i_1,i_2,\dots i_k\}\subseteq\{1,2,\dots,n\}$) can be represented by the tuples $(z_{1},z_{2},\dots, z_{n})$ where 
$$
z_i=\begin{cases}
t & \text{if $i\in \{i_1,i_2,\dots , i_k\}$ }\\
0 & \text{otherwise}.
\end{cases}
$$
The blow-up $L^B$ of $L\cong \mathbf{2}^3$ is shown in Figure \ref{figure3}.

\textbf{Throughout this paper, $L^B$ denotes the blow-up of a Boolean lattice $L\cong \mathbf{2}^n $ with $n\geq 3$ and hence  $|L^B|\geq 8$. }

\begin{remark}\label{pcbool}
Note that if $a, b \in L$ ($a\not=b$), where $L$ is Boolean and $C_a$ and $C_b$ are the corresponding chains in $L^B$, then $a\wedge b =x \wedge y $ and $a\vee b =x \vee y $ in $L^B$  for any element $x$ on the chain $C_a$ and any element $y$ on the chain $C_b$. Hence, in particular, if $x\in L^B$ and $x^*$ be the pseudocomplement of $x$ in $L^B$, then in $L^B$, we have $x \vee x^*=1$ and $x \wedge x^* =0$. Note that $x^*$ need not be the unique complement of $x$ in $L^B$, whereas $x^*$ is the unique complement of $x$ in $L$. Also, in $L^B$, the pseudocomplement of atom $q_i$ is the dual atom of $L^B$, denoted by $q_i^*$ and the pseudocomplement of dual atom $q_i^*$ is the largest element in the chain of $[q_i]$.
\end{remark}
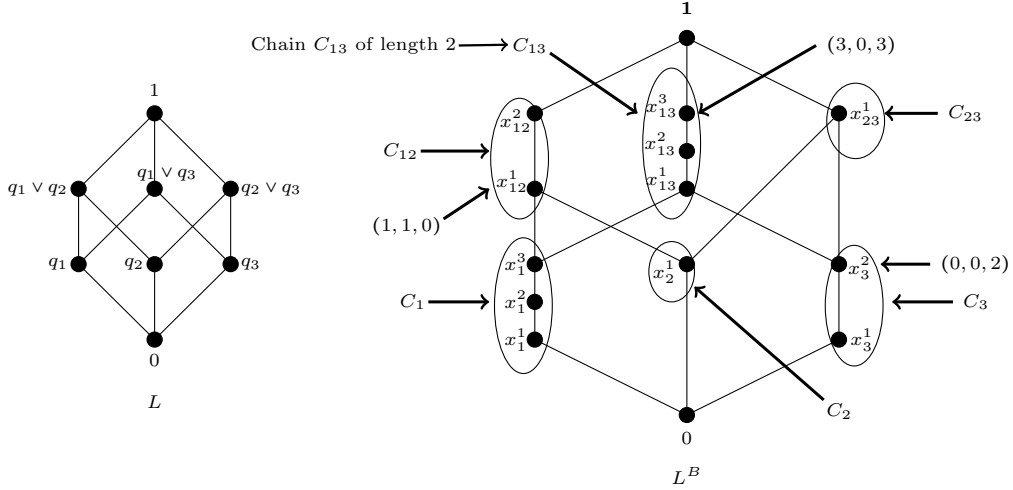
\begin{figure}[h]
\begin{center}
\begin{tikzpicture}	[scale=1]		

\begin{scope}[shift={(-2,1)}]
	\draw (0,0) -- (-1,1); \draw (0,0) -- (0,1); \draw (0,0) -- (1,1);
	\draw (-1,1) -- (-1,2); \draw (-1,1) -- (0,2); \draw (0,1) --
	(-1,2); \draw (0,1) -- (1,2); \draw (1,1) -- (1,2); \draw (1,1) --
	(0,2);
	\draw (-1,2) -- (0,3); \draw (0,2) -- (0,3); \draw (1,2) -- (0,3);
	\draw [fill=black] (0,0) circle (.1); \draw [fill=black] (0,1)
	circle (.1);\draw [fill=black] (0,2) circle (.1);\draw
	[fill=black] (0,3) circle (.1);
	\draw [fill=black] (-1,1) circle (.1); \draw [fill=black] (-1,2)
	circle (.1);
	\draw [fill=black] (1,1) circle (.1);\draw [fill=black] (1,2)
	circle (.1);
	\node [below] at (0,-.05) {\tiny$0$};\node [left] at (-1,1)
	{\tiny$q_1$};\node [right] at (1,1) {\tiny$q_3$}; \node [left] at (0,1)
	{\tiny$q_2$};
	\node [left] at (-1,2) {\tiny$q_1\vee q_2$}; \node [left] at (0.7,2.2) {\tiny$q_1\vee q_3$};
	\node [right] at (1,2) {\tiny$q_2\vee q_3$};
	\node [above] at (0,3.1) {\tiny$1$};
	\node[below] at (0,-0.6) {\tiny$L$};
	
\end{scope}

\begin{scope}[shift={(5,0)}]
	
	\draw [fill=black] (0,0) circle (.1);
	\draw [fill=black] (-2,1) circle (.1);
	
	\draw [fill=black] (-2,2) circle (.1);
	
	\draw [fill=black] (-2,3) circle (.1);
	
	\draw [fill=black] (-2,4) circle (.1);

	\draw [fill=black] (0,2) circle (.1);
	\draw [fill=black] (0,4) circle (.1);
	\draw [fill=black] (-2,1.5) circle (.1);
	\draw [fill=black] (0,3.5) circle (.1);
	\draw [fill=black] (0,3) circle (.1);
	
	\draw [fill=black] (2,2) circle (.1);
	\draw [fill=black] (2,4) circle (.1);
	\draw [fill=black] (2,1) circle (.1);
	\draw [fill=black] (0,5) circle (.1);
	
	\draw (0,0) -- (-2,1) --  (-2,2) -- (-2,3) -- (-2,4)--(0,5)-- (0,4) --(0,3) --(-2,2);
	\draw (0,0) -- (0,2) --(-2,3);
	\draw (0,0) -- (2,1)--(2,2) --(2,4) --(0,5);
	\draw (0,2) -- (2,4);
	\draw (0,3) -- (2,2);

	\draw (-2.15,1.45) ellipse (0.38 and 0.9); \draw (-2.2,3.4) ellipse (0.38 and 0.8); \draw (2.2,1.45) ellipse (0.38 and 0.8); \draw (-0.2,1.9) ellipse (0.3 and 0.4); \draw (-0.2,3.6) ellipse (0.38 and 1);\draw (2.22,3.9) ellipse (0.38 and 0.5);
	
	\draw [line width=0.40mm,->] (-3.4,1.5)--(-2.6,1.5);
	\node [left] at (-3.3,1.5) {\tiny $C_{1}$};
	
	\draw [line width=0.40mm,->] (-3.5,3.5)--(-2.6,3.5);
	\node [left] at (-3.4,3.5) {\tiny $C_{12}$};
	
	\draw [line width=0.40mm,<-] (2.7,1.5)--(3.5,1.5);
	\node [right] at (3.5,1.5) {\tiny $C_{3}$};
	
	\draw [line width=0.40mm,<-] (2.6,4)--(3.3,4);
	\node [right] at (3.3,4) { \tiny$C_{23}$};
	
	\draw [line width=0.40mm,<-] (-0.65,4)--(-1.8,4.8);
	\node [left] at (-1.7,4.9) {\tiny$C_{13}$};
	
	\draw [line width=0.30mm,->] (-3,4.9)--(-2.35,4.91);
	\node [left] at (-2.9,4.9) {\tiny Chain $C_{13}$ of length $2$};

	\draw [line width=0.40mm,<-] (0.1,1.7)--(1.8,0.2);
	\node [right] at (1.7,0.05) {\tiny $C_{2}$};
	
	\draw [line width=0.40mm,<-] (2.55,2)--(3.2,2);
	\node [right] at (3.2,2) { \tiny$(0,0,2)$};

	\draw [line width=0.40mm,<-] (0.15,4)--(1.7,4.9);
	\node [right] at (1.7,4.9) { \tiny$(3,0,3)$};
	
	\draw [line width=0.40mm,->] (-3.2,2.6)--(-2.6,3);
	\node [left] at (-3.1,2.5) { \tiny$(1,1,0)$};

	\node at (-2.25,1) {\tiny$ x_1 ^{1}$};
	\node at (-2.25,1.5) {\tiny$ x_1 ^{2}$};
	\node at (-2.25,2) {\tiny$ x_1 ^{3}$};
	\node at (-0.3,1.9) {\tiny$ x_2 ^{1}$};
	\node at (2.27,1.95) {\tiny$ x_3 ^{2}$};
	\node at (-2.3,3.1) {\tiny$ x_{12} ^{1}$};
	\node at (-2.25,3.9) {\tiny$ x_{12} ^{2}$};
	\node at (-0.33,4.1) {\tiny$ x_{13} ^{3}$};
	\node at (-0.35,3.6) {\tiny$ x_{13} ^{2}$};
	\node at (-0.33,3.13) {\tiny$ x_{13} ^{1}$};
	\node at (2.35,4) {\tiny$ x_{23} ^{1}$};
	\node at (2.3,1) {\tiny$ x_3 ^{1}$};
	
	\node at (0,-0.8) {\tiny$ L^B$};

	\node at (0,-0.3) {\tiny 0};
	\node at (0,5.4) {\tiny $\mathbf{1}$};
\end{scope}

\end{tikzpicture}
\caption{Boolean lattice $L\cong \mathbf{2}^3$ and its blow-up $L^B$ }\label{figure3}
\end{center}
\end{figure}

In particular, in $L^B$, we observe that, $[x_{i_1i_2 \dots i_k}^{1}]= [x_{i_1i_2 \dots i_k}^{2}]= \dots = [x_{i_1i_2 \dots i_k}^{n_j}]$, where $\{i_1,i_2, \dots ,i_k\}\subseteq \{1,2,\dots,n\}\}$. Thus, the elements on the chain $C_{12 \dots k}$ have the same equivalence classes.
The following result is due to Gadge et al. \cite{pgg}.

\begin{lemma}[{Gadge et al. \cite[Lemma 3.8]{pgg}}]\label{ll'}
Let $L$ be a pseudocomplemented lattice and $L'$ be a poset obtained from $L$ by replacing an element of $L$ with a bounded chain. Then  $L'$ is pseudocomplemented.
\end{lemma}

\begin{corollary}\label{pseudo}
	Let $L^B$ be a blow-up of a Boolean lattice $L\cong \mathbf{2}^n$. Then the following statements  hold: \begin{enumerate}
		\item $L^B$ and its dual lattice $(L^B)^{\partial}$ both are pseudocomplemented.
		\item $[L^B]\cong [(L^B)^\partial]\cong L\cong\mathbf{2}^n$.
		\item Let $a,b\in L^B$. Then $a^\perp=b^\perp$ if and only if $a^*=b^*$, where $a^*,b^*$ denotes the pseudocomplement of $a$ and $b$ in $L^B$ respectively.
	\end{enumerate}
\end{corollary}

\begin{proof}
	\begin{enumerate}
		\item The proof follows from Lemma \ref{ll'}.
		\item It is well known that every pseudocomplemented lattice is $0$-distributive. By (1), $L^B$ and $(L^B)^\partial$ both are pseudocomplemented. Hence by Theorem \ref{paka}, $[L^B]$ and $[(L^B)^\partial]$ are Boolean.
		\item Suppose that $a^\perp=b^\perp$. Let $a^*$ be the pseudocomplement of $a$  in $L^B$. Then $a\wedge a^*=0$. This implies that $a^*\in a^\perp=b^\perp$. Thus $a^* \wedge b=0$. This shows that $a^*\leq b^*$, where $b^*$ is the pseudocomplement of $b$. Similarly, we can show that $b^*\leq a^*$. Hence $a^*=b^*$.
		\par Conversely, assume that $a^*=b^*$. Let $x \in a^\perp$. Then $x \wedge a=0$. Hence $x \leq a^*=b^*$. This further gives $x \wedge b=0$. Hence, $x \in b^\perp$. Thus $a^\perp \subseteq b^\perp$. Similarly, we can prove that  $b^\perp \subseteq a^\perp$. This proves that $a^\perp =b^\perp$.
\end{enumerate} \end{proof}

Now, we prove that the zero-divisor graph of a 0-distributive lattice with $n$ atoms can be realized as the zero-divisor graph of a blow-up of a Boolean lattice $L=\mathbf{2}^n$.

\begin{theorem}\label{0-disblowup}
	Let $L'$ be a finite 0-distributive lattice with $n$ atoms. Let $L^B$ be the blow-up of the Boolean lattice $L=\mathbf{2}^n$. Then $G(L')=G(L^B)$.
\end{theorem}
\begin{proof}
	By Theorem \ref{paka}, we have $[L']$ is Boolean and $[L']\cong \mathbf{2}^n$, as $L'$ has $n$ atoms. Note that $L \cong [L']$. Let $q_i$ $(1 \leq i \leq n)$ be all $n$ atoms of $L'$. Let $|[q_i]|=m_i$ for every $i$. To construct $L^B$ from $L$, replace each atom, say $p_i$, of $L$, by the chain of length $m_i$. Further, let $[x]$ be an element of $[L']$. Since $L \cong [L']$, we can assume that $x$ is an element of $L$, which is an image of $[x]\in [L']$. Now, replace $x$ in $L$ with the chain of length $|[x]|$. It is easy to observe that if $d$ is a dense element of $L'$, then clearly $\{d,1\}\subseteq [d]$. In this case, the elements in $[d]$ will be dense in $L^B$ too. This gives that $Z(L')=Z(L^B)$. Hence, the vertex sets of  $G(L')$ and $G(L^B)$ are the same. Further, one can see that $x \wedge y=0$ in $L'$ if and only if $[x] \wedge [y]=[0]=[x\wedge y]$ in $[L']$.  So $x$ and $y$ can be viewed as elements of $L$ and hence $x \wedge y=0$ in $L$. Let $C_x$ and $C_y$ be the chains in $L^B$ obtained by replacing $x$ and $y$ in $L$. By Remark \ref{pcbool}, we have $a \wedge b =0$ for every $a \in C_x$ and every $b \in C_y$. Hence, in particular, in $L^B$ also, $x \wedge y=0$. The converse follows on similar lines. This proves that $G(L')=G(L^B)$.
\end{proof}

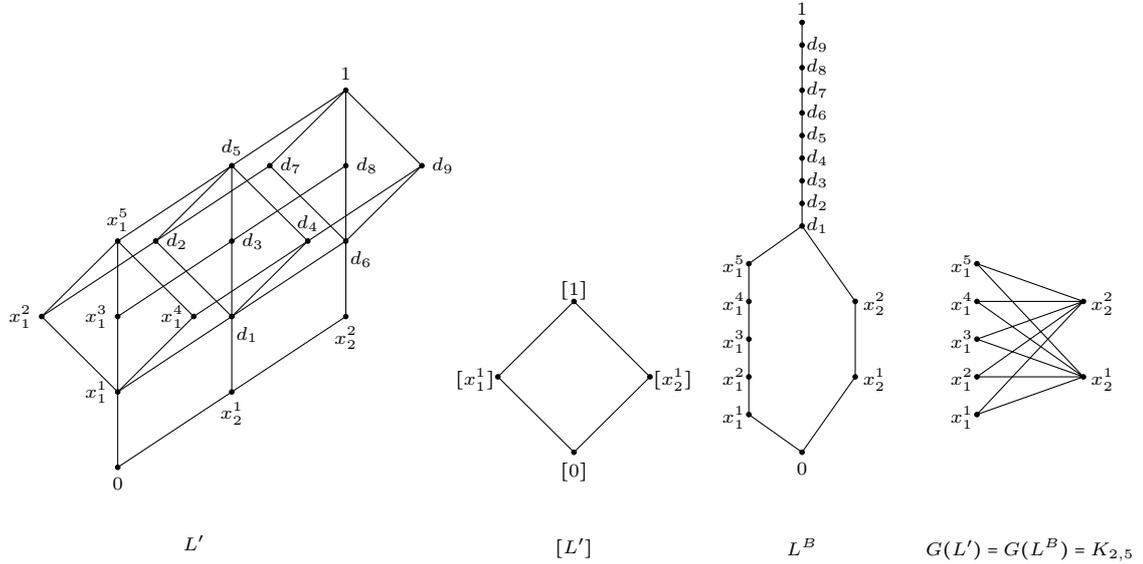
\begin{figure}[h]
	\begin{center}
		
		\begin{tikzpicture}[scale=1]
			
			\begin{scope}[shift={(-8,0)}]
				\draw (0,0) -- (0,1); \draw (0,1) -- (-1,2); \draw (0,1) -- (0,2);
				\draw (0,1) -- (1,2); \draw (-1,2) -- (0,3); \draw (0,2) -- (0,3);
				\draw (1,2) -- (0,3);
				
				\draw (1.5,1) -- (1.5,2); \draw (1.5,2) -- (.5,3); \draw (1.5,2)
				-- (1.5,3); \draw (1.5,2) -- (2.5,3); \draw (.5,3) -- (1.5,4);
				\draw (1.5,3) -- (1.5,4); \draw (2.5,3) -- (1.5,4);

				\draw (3,2) -- (3,3); \draw (3,3) -- (2,4); \draw (3,3) -- (3,4);
				\draw (3,3) -- (4,4); \draw (2,4) -- (3,5); \draw (3,4) -- (3,5);
				\draw (4,4) -- (3,5);
				
				\draw (0,0) -- (1.5,1); \draw (1.5,1) -- (3,2);
				
				\draw (0,1) -- (1.5,2); \draw (1.5,2) -- (3,3);
				
				\draw (-1,2) -- (.5,3); \draw (0.5,3) -- (2,4);
				
				\draw (0,2) -- (1.5,3); \draw (1.5,3) -- (3,4);
				
				\draw (1,2) -- (2.5,3); \draw (2.5,3) -- (4,4);
				
				\draw (0,3) -- (1.5,4); \draw (1.5,4) -- (3,5);
				
				\draw [fill=black] (0,0) circle (.03); \draw [fill=black] (0,1)
				circle (.03);\draw [fill=black] (0,2) circle (.03);\draw
				[fill=black] (0,3) circle (.03);
				
				\draw [fill=black] (-1,2) circle (.03);\draw [fill=black] (1,2)
				circle (.03);
				
				\draw [fill=black] (1.5,1) circle (.03);\draw [fill=black] (1.5,2)
				circle (.03);\draw [fill=black] (1.5,3) circle (.03);\draw
				[fill=black] (1.5,4) circle (.03);
				
				\draw [fill=black] (.5,3) circle (.03);\draw [fill=black] (2.5,3)
				circle (.03);
				
				\draw [fill=black] (3,2) circle (.03);\draw [fill=black] (3,3)
				circle (.03);\draw [fill=black] (3,4) circle (.03);\draw
				[fill=black] (3,5) circle (.03);
				
				\draw [fill=black] (2,4) circle (.03);\draw [fill=black] (4,4)
				circle (.03);
				
				\node [below] at (0,0) {\tiny $0$}; \node [left] at (0,1)
				{\tiny $x_1^{1}$}; \node [left] at (-1,2) {\tiny $x_1^{2}$}; \node
				[left] at (0,2) {\tiny $x_1^{3}$}; \node [left] at (1,2) {\tiny
					$x_1^{4}$}; \node [above] at (0,3) {\tiny $x_1^{5}$};
				
				\node [below] at (1.5,1) { \tiny $x_2^{1}$}; \node [below] at
				(1.7,2) {\tiny $d_1$}; \node [right] at (1.5,3) {\tiny $d_3$};
				\node [above] at (1.5,4) {\tiny $d_5$}; \node [right] at (0.5,3)
				{\tiny $d_2$}; \node [above] at (2.5,3) {\tiny $d_4$};
				
				\node [below] at (3,2) {\tiny $x_2^{2}$}; \node [below] at (3.2,3)
				{\tiny $d_6$}; \node [right] at (3,4) {\tiny $d_8$}; \node
				[above] at (3,5) {\tiny $1$}; \node [right] at (2,4) {\tiny
					$d_7$}; \node [right] at (4,4) {\tiny $d_9$};
				
				\node [] at (1,-1) {\tiny$L'$};
			\end{scope}
			
			\begin{scope}[shift={(-2,-0.8)}]
				\draw [fill=black] (0,1)
				circle (.03);
				\draw	[fill=black] (0,3) circle (.03);
				
				\draw [fill=black] (-1,2) circle (.03);
				\draw [fill=black] (1,2)
				circle (.03);
				
				\draw (0,1)--(-1,2)--(0,3)--(1,2)--(0,1);
				
				\node [below] at (1.3,2.25) {\tiny $[x_2^{1}]$};
				
				\node [below] at (-1.3,2.25) {\tiny $[x_1^{1}]$};
				
				\node [below] at (0,3.4) {\tiny $[1]$};
				
				\node [below] at (0,1) {\tiny $[0]$};
				
				\node [below] at (0,0) {\tiny $[L']$};
			\end{scope}
			
			\begin{scope}[shift={(1,-0.8)}]
				\draw [fill=black] (0,1)
				circle (.03);
			
				\draw	[fill=black] (0.7,3) circle (.03);
				
				\draw [fill=black] (-0.7,3) circle (.03);
				\draw [fill=black] (-0.7,3.5) circle (.03);
				\draw [fill=black] (-0.7,1.5) circle (.03);
				\draw [fill=black] (-0.7,2) circle (.03);
				\draw [fill=black] (-0.7,2.5) circle (.03);
				\draw [fill=black] (0.7,2)
				circle (.03);
				\draw [fill=black] (0,4) circle (.03);
					\draw [fill=black] (0,4.3) circle (.03);
					\draw [fill=black] (0,4.6) circle (.03);
					\draw [fill=black] (0,4.9) circle (.03);
					\draw [fill=black] (0,5.2) circle (.03);
					\draw [fill=black] (0,5.5) circle (.03);
					\draw [fill=black] (0,5.8) circle (.03);
					\draw [fill=black] (0,6.1) circle (.03);
					\draw [fill=black] (0,6.4) circle (.03);
					\draw [fill=black] (0,6.7) circle (.03);
				
				\draw (0,1)--(-0.7,1.5)--(-0.7,2)--(-0.7,2.5)--(-0.7,3)--(-0.7,3.5)--(0,4)--(0.7,3)--(0.7,2)--(0,1);
				\draw (0,4)--(0,4.3)--(0,4.6)--(0,4.9)--(0,5.2)--(0,5.5)--(0,5.8)--(0,6.1)--(0,6.4)--(0,6.7);
				
				\node [below] at (0.95,2.25) {\tiny $x_2^{1}$};
				
				\node [below] at (0.95,3.25) {\tiny $x_2^{2}$};
				
				\node [below] at (-0.9,2.25) {\tiny $x_1^{2}$};
				\node [below] at (-0.9,1.75) {\tiny $x_1^{1}$};
				\node [below] at (-0.9,2.75) {\tiny $x_1^{3}$};
				\node [below] at (-0.9,3.25) {\tiny $x_1^{4}$};
				\node [below] at (-0.9,3.75) {\tiny $x_1^{5}$};
				
				\node [below] at (0,7.1) {\tiny $1$};
				\node [below] at (0.2,4.25) {\tiny $d_1$};
				\node [below] at (0.2,4.55) {\tiny $d_2$};
				\node [below] at (0.2,4.85) {\tiny $d_3$};
				\node [below] at (0.2,5.15) {\tiny $d_4$};
				\node [below] at (0.2,5.45) {\tiny $d_5$};
				\node [below] at (0.2,5.75) {\tiny $d_6$};
				\node [below] at (0.2,6.05) {\tiny $d_7$};
				\node [below] at (0.2,6.35) {\tiny $d_8$};
				\node [below] at (0.2,6.65) {\tiny $d_9$};
				
				\node [below] at (0,1) {\tiny $0$};
				
				\node [below] at (0,0) {\tiny $ L^B$};
			\end{scope}
			
				\begin{scope}[shift={(4,-0.8)}]

				\draw	[fill=black] (0.7,3) circle (.03);
				\draw [fill=black] (-0.7,3) circle (.03);
				\draw [fill=black] (-0.7,3.5) circle (.03);
				\draw [fill=black] (-0.7,1.5) circle (.03);
				\draw [fill=black] (-0.7,2) circle (.03);
				\draw [fill=black] (-0.7,2.5) circle (.03);
				\draw [fill=black] (0.7,2)
				circle (.03);

				\draw (-0.7,1.5)--(0.7,2)--(-0.7,2)--(0.7,3)--(-0.7,2.5)--(0.7,2)--(-0.7,3)--(0.7,3)--(-0.7,3.5)--(0.7,2);
				 \draw (-0.7,1.5)--(0.7,3);

				\node [below] at (0.95,2.25) {\tiny $x_2^{1}$};
				
				\node [below] at (0.95,3.25) {\tiny $x_2^{2}$};
				
				\node [below] at (-0.9,2.25) {\tiny $x_1^{2}$};
				\node [below] at (-0.9,1.75) {\tiny $x_1^{1}$};
				\node [below] at (-0.9,2.75) {\tiny $x_1^{3}$};
				\node [below] at (-0.9,3.25) {\tiny $x_1^{4}$};
				\node [below] at (-0.9,3.75) {\tiny $x_1^{5}$};

				\node [below] at (0,0) {\tiny $G(L')=G(L^B)=K_{2,5}$};
			\end{scope}
			
		\end{tikzpicture}
		\caption{Illustration of Theorem \ref{0-disblowup} }\label{p1.p2}
	\end{center}
	
\end{figure}

\begin{remark}\label{rem4.6}
	Since every pseudocomplemented lattice is 0-distributive and the converse is true if a lattice is finite, we have, by Corollary \ref{pseudo}, $L^B$ is a $0$-distributive lattice and hence for every  $x,y \in L^B$, $x\vee y \in Z^*(L^B)$  if and only if $x^\perp \cap y^\perp \neq\{0\}$. In fact, in a 0-distributive lattice, $(x\vee y)^\perp=x^\perp\cap y^\perp$.
\end{remark}

The following result essentially follows from Alizadeh et al. in \cite{aliz}. 
\begin{lemma}[{M. Alizadeh et al. \cite[Theorem 3.3]{aliz}}] \label{alizadeh}\label{distance}
	Let $L^B$ be a blow-up of a Boolean lattice $L\cong \mathbf{2}^n$. Then for $x,y\in V(G(L^B))$ the following assertions hold:
	\begin{enumerate} 
		\item  $d(x,y)=1$ if and only if $y\in (x^{*}]$.
		\item  $d(x,y)=2$ if and only if $y\notin (x^*]$ and $y^*\notin (x^{**}]$.
		
		\item  $d(x,y)=3$ if and only if $y\notin (x^*]$ and $y^*\in (x^{**}]$.
	\end{enumerate}
\end{lemma}

\begin{lemma}\label{mutually}
	Let $L^B$ be the blow-up of a Boolean lattice $L\cong \mathbf{2}^n$. If $x$ and $y$ are comparable or $x\wedge y =0$, then $x$ and $y$ are not mutually maximally distant for every pair $a,b \in V(G(L^B))$ with  $x\in [a]$, $y\in [b]$ and $[a]\neq [b]$.
\end{lemma}
\begin{proof}
	Let $x\in [a]$ and $y\in [b]$ for $a,b \in V(G(L^B))$. If $x\wedge y=0$, then $d(x,y)=1$. Now, by the assumption $|L^B|\geq 8$, and by Remark \ref{diam}, we can easily find $z\in N(x)$ such that $d(z,y)= 2$. This implies that $x$ and $y$ are not mutually maximally distant.
	\par Let $x$ and $y$ be comparable. Without loss of generality, assume that $x\leq y$. Then  $x\vee y \neq 1$ and $x\wedge y \neq 0$, as $x,y \in V(G(L^B))$. By Lemma \ref{alizadeh}(2), $d(x,y)=2$. 
	
	To prove that $x$ and $y$ are not mutually maximally distant, we take $z=x^*$. Then we have $z\in N(x)$. We show that $d(y,z)=3$.
	
	By Remark \ref{pcbool}, $x\vee z=1$ and $x\leq y$, we have $y\vee z=1$. Also, observe that $y\wedge z\neq 0$, that is, $y\notin (z^*]$ and hence $d(y,z)\neq 1$.
	
	For this, if $y\wedge z= 0$, i.e., $y\wedge x^*=0$, then $x^*\leq y^*$, as $L^B$ is pseudocomplemented. Further, as $x\leq y$, we have $y^*\leq x^*$. Thus $x^*=y^*$. By Corollary \ref{pseudo} (3), $x^\perp =y^\perp$, and hence $[x]=[y]$, a contradiction, as $x\in [a]$ and $y\in [b]$ with $[a]\neq [b]$.  Hence $y\wedge z\not = 0$.
	
	We claim that $d(y,z)\neq 2$. For this, if $d(y,z)=2$, there is a path $y-t-z$ in $G(L^B)$. Therefore $y\wedge t=0=z\wedge t$. As $L^B$ is 0-distributive, we have $0=t\wedge (y\vee z)=t\wedge 1=t=0$, a contradiction to $t\in V (G(L^B))$.

	So by part (3) of Lemma \ref{alizadeh}, $d(y,z)=3$, as $y\vee z =1$ gives $y^* \wedge z^*=0$, that is, $y^*\in (z^{**}]$. This means that $x$ and $y$ are not mutually maximally distant.
\end{proof}
\begin{definition}\label{GLB}

	We associate a graph  $G(L^B)^{**}$ with the lattice $L^B$ whose  vertex set is   $V(G(L^B)^{**})=V(G(L^B))$ and  two distinct vertices $x$ and $y$ are adjacent in $G(L^B)^{**}$ if and only if either $([x]=[y])$ or $(a\wedge b \neq 0$, $a\nleq b$ and  $b\nleq a$ for every $a\in [x]$ and for every $b\in [y])$.  Equivalently, we can prove that $x$ and $y$ are adjacent in $G(L^B)^{**}$ if and only if either $([x]=[y])$ or ($[x] \wedge [y]\neq [0]$ with $[x] || [y]$). 
	
	Also, let $G(L^B)^{*}=G(L^B)$, if $G(L^B)$ is complete, otherwise, $G(L^B)^{*}$ is extracted from $G(L^B)^{**}$ after deleting all isolated vertices. Note that if $x$ and $y$ are vertices of $G(L^B)$ with $[x]=[y]$, then $N(x)=N(y)$.
	
\end{definition}

\begin{lemma}	\label{isoatom}
	Suppose $L^B$ be the blow-up of a Boolean lattice $L\cong \mathbf{2}^n$ ($n \geq 3$) by replacing the elements of $L$ by chains of finite length except for the atoms $q_i$ ($1 \leq I \leq n$) of $L$. Then  $x$ is an isolated vertex in $G(L^B)^{**}$ if and only if $x\in Atoms(L^B)$.	
\end{lemma}

\begin{proof}
	Suppose on the contrary that $x$ is not an isolated vertex in $G(L^B)^{**}$ for some $x\in Atoms(L^B)$. Then $x$ is adjacent to some $y$ in $G(L^B)^{**}$. Then by the definition of adjacency in $G(L^B)^{**}$, we have either ($x \wedge y \neq 0$, $x \not\leq y$ and $y \not\leq x$) or ($[x]=[y]$). If $x \wedge y \neq 0$,  then we have $x \leq y$, as $x\in Atoms(L^B)$, a contradiction to $x \not\leq y$. 	Hence $[x]=[y]$, that is, $x^\perp =y^\perp$. 
	
	Since $x\in Atoms(L^B)$, we have $x$ is incomparable with $y$ or $x \leq y$. If $x$ is incomparable with $y$, then $x \in y^\perp=x^\perp$, a contradiction. 	Hence $x \leq y$.

	If $x < y$ and if there is some atom $q< y$, then $x^\perp \neq y^\perp$, a contradiction to the fact that $q\in x^\perp$ and $q\notin y^\perp$. Hence, there is no atom $q < y$. However, in this case, $x<y$ forms a chain above the atom $x$. This contradicts the fact that, in $L^B$, atoms are not replaced by chains. Thus $x=y$. 
	This is again a contradiction to $x \neq y$. Hence, $[x]\neq [y]$. 
	
	Hence $x$ is an isolated vertex in $G(L^B)^{**}$.		
	
	Conversely, assume that  $x$ is an isolated vertex in $G(L^B)^{**}$. 	Suppose $x\notin Atoms(L^B)$.
	
	Let $C_{12\dots (i-1)(i+1) \cdots n}$ be a chain in $L^B$ obtained by replacing a dual atom $d_i= \bigvee\limits_{j=1, i\neq j}^{j=n} q_j$, where $q_j$ are atoms of $L^B$. 
	
	Let $S=\bigl\{x^{t}_{12\dots (i-1)(i+1) \cdots n}~| \text{ for some }x^{t}_{12\dots (i-1)(i+1) \cdots n} \in C_{12\dots (i-1)(i+1) \cdots n},  \forall i, 1 \leq i \leq n \bigr\}$.
	
	Clearly, for given $x\in V(G(L^B)^{**}) \setminus Atoms(L^B)$, there is a $y \in L^B$ such that $x \wedge y \neq 0, ~ x \not\leq y$ and $y \not\leq x$. Note that such $y$ exists, as $n\geq 3$ and $y$ be an element of the set $S$, in particular, which is above the pseudocomplement of $x^*$. Thus, $x$ is adjacent to $y$, a contradiction to the fact that $x$ is an isolated vertex in $G(L^B)^{**}$. Hence $x\in Atoms(L^B)$.
\end{proof}

The following result gives the structure of $G(L^B)^{**}$.

\begin{lemma}\label{blowl2}
	Suppose $L^B$ is a blow-up of a Boolean lattice $L\cong \mathbf{2}^n$ with $n\geq 3$. Then $G(L^B)^{**}= H+ K_{|[q_{1}]|}+K_{|[q_{2}]|}+ \dots +K_{|[q_{n}]|}$, where $H$ is a connected graph and $K_{|[q_i]|}$ be the complete graph on $|[q_i]|$ vertices.
\end{lemma}
\begin{proof}
	Let $q_i$, $ 1 \leq i \leq n$ be the all atoms of $L^B$. Note that $[q_i]\neq [q_j]$ for $i \neq j$.	
	Let $A_1=\{x\in V(G(L^B)) \mid x\in [q_1]\}$, $A_2=\{x\in V(G(L^B)) \mid x\in [q_2]\}$,\dots, $A_n=\{x\in V(G(L^B)) \mid x\in [q_n]\}$ and $A=\bigcupdot A_i$, $1\leq i\leq n$. 
	
	\par Now, we partition the vertex set of $G(L^B)^{**}$ as $V(G(L^B)^{**})=V(G(L^B))\setminus A$. If $x,y\in A_i$, then we have $[x]=[y]$ and hence $x$ is adjacent to $y$ in $G(L^B)^{**}$. This implies that $G(L^B)^{**}[A_i]$ is a complete graph, for every $1\leq i\leq n$.

	Suppose that $x\in A_i$ and $y \in A_j$ for $i\neq j$. Then $[x]\neq [y]$. This implies that $a\wedge b =0$ for every $a\in [x]$ and $b\in [y]$, and so $x$ is not adjacent to $y$ in $G(L^B)^{**}$. Suppose $x\in A_i$ and $y\in V(G(L^B))\setminus A$. Hence  $[x]\neq [y]$. Also, for every $a\in [x]$, we can easily observe  that either $a\wedge b =0$ or $a\leq b $ for every $b\in [y]$. This means that $x$ and $y$ are not adjacent in $G(L^B)^{**}$.

	Consider a graph $H$ with the vertex set as $V(H)= V(G(L^B))\setminus A$. We prove that $H$ is a connected graph. 
	For this, consider a set defined in Lemma \ref{isoatom} as
	
	 $S=\bigl\{x^{t}_{12\dots (i-1)(i+1) \cdots n}~| \text{ for some }x^{t}_{12\dots (i-1)(i+1) \cdots n} \in C_{12\dots (i-1)(i+1) \cdots n},  \forall i, 1 \leq i \leq n \bigr\}$.

	Now, for given $x\in V(G(L^B)^{**}) \setminus Atoms(L^B)$, there is a $y \in L^B$ such that $x \wedge y \neq 0, ~ x \not\leq y$ and $y \not\leq x$. Note that such $y$ exists, as $n\geq 3$ and $y$ be an element of the set $S$, in particular, which is above the pseudocomplement of $x^*$.

	
	Thus, $H$ is a connected graph. Hence $G(L^B)^{**}= H+ K_{|[q_{1}]|}+K_{|[q_{2}]|}+ \dots +K_{|[q_{n}]|}$. 
\end{proof}	

Consider the lattice $L$ and $L^B$ shown in Figure \ref{figure3}. The graph $G(L^B)^{**}$ is shown in  Figure \ref{figure4}. This illustrates  Lemma \ref{blowl2}.

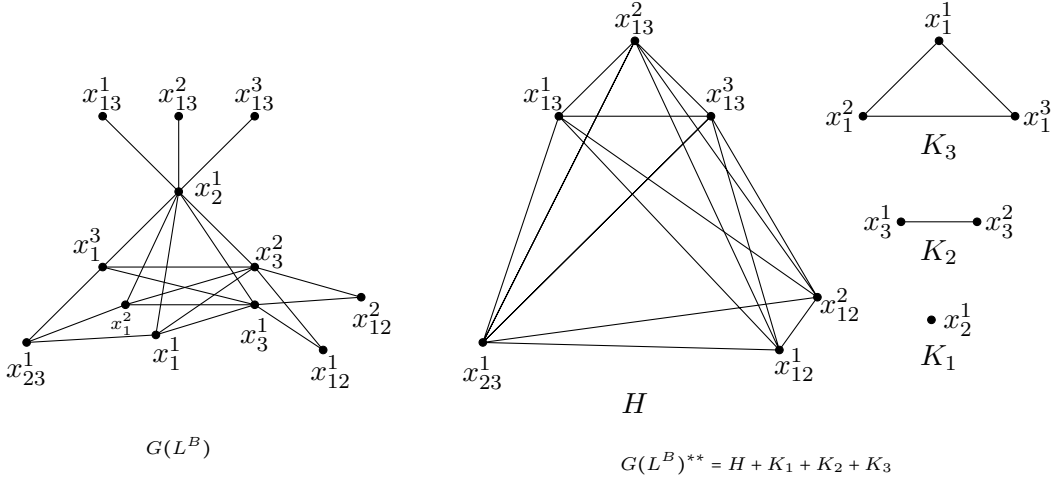
\begin{figure}[h]
	\begin{center}
		\begin{tikzpicture}	[scale=1]

			\begin{scope}[shift={(-4,0)}]
				\draw [fill=black] (-2,1) circle (.05);
				\draw [fill=black] (-1,2) circle (.05);
				
				\draw [fill=black] (1,2) circle (.05);

				\draw [fill=black] (-0.7,1.5) circle (.05);
				
				\draw [fill=black] (-0.3,1.1) circle (.05);
				
				\draw [fill=black] (1.9,0.9) circle (.05);
				\draw [fill=black] (2.4,1.6) circle (.05);
				\draw [fill=black] (-1,4) circle (.05);
				\draw [fill=black] (1,4) circle (.05);
				\draw [fill=black] (0,4) circle (.05);
				\draw [fill=black] (0,3) circle (.05);
				\draw [fill=black] (1,1.5) circle (.05);
				
				\node at (0.4,3.051) {$ x_2 ^{1}$};
				\node at (1.2,2.2) {\bf$ x_3 ^{2}$};
				\node at (-0.15,0.85) {$ x_1 ^{1}$};
				\node at (-0.75,1.25) {\tiny$ x_1 ^{2}$};
				\node at (-1.2,2.3) {$ x_1 ^{3}$};
				\node at (2.0,0.61) {$ x_{12} ^{1}$};
				\node at (2.55,1.32) {\bf$ x_{12} ^{2}$};
				\node at (-1,4.3) {\bf$ x_{13} ^{1}$};
				\node at (1,4.3) {\bf$ x_{13} ^{3}$};
				\node at (-2,0.65) {\bf$ x_{23} ^{1}$};
				\node at (0,4.3) {\bf$ x_{13} ^{2}$};
				\node at (1,1.1) {\bf$ x_3 ^{1}$};
				
				\draw (-2,1)--(-0.3,1.1) -- (1,2);
				\draw (-0.3,1.1) -- (1,1.5)--(2.4,1.6);
				\draw (1,1.5)--(1.9,0.9);
				\draw (1,1.5)--(0,3);
				\draw (-0.7,1.5)--(-2,1);
				\draw (-0.7,1.5) -- (1,1.5);
				\draw (1,1.5) -- (-1,2);
				\draw (-0.7,1.5) -- (1,2);
				\draw (-0.7,1.5) -- (0,3)--(-0.3,1.1);
				\draw (0,3)--(1,4);
				\draw (0,3)--(0,4);
				\draw (0,3)--(-1,4);
				\draw (0,3)-- (-1,2)-- (1,2)--(0,3);
				\draw (-2,1)--(-1,2);
				\draw (1,2)--(1.9,0.9);
				\draw (1,2)--(2.4,1.6);
				\node at (0,-0.4) {\tiny $ G(L^B)$};

			\end{scope}
			\begin{scope}[shift={(2,0)}]
				\draw [fill=black] (-2,1) circle (.05);

				\draw [fill=black] (1.9,0.9) circle (.05);
				\draw [fill=black] (2.4,1.6) circle (.05);
				\draw [fill=black] (-1,4) circle (.05);
				\draw [fill=black] (1,4) circle (.05);
				\draw [fill=black] (0,5) circle (.05);
				
				\draw [fill=black] (4,5) circle (.05);
				\draw [fill=black] (3,4) circle (.05);
				\draw [fill=black] (5,4) circle (.05);
				\draw [fill=black] (4.5,2.6) circle (.05);
				\draw [fill=black] (3.5,2.6) circle (.05);
				\draw [fill=black] (3.9,1.3) circle (.05);
				
				\draw (4,5)--(3,4)--(5,4)--(4,5);
				\draw (4.5,2.6)--(3.5,2.6);
				
				\node at (4,5.3) {\bf$ x_{1} ^{1}$};
				\node at (2.7,4) {\bf$ x_{1} ^{2}$};
				\node at (5.3,4) {\bf$ x_{1} ^{3}$};
				\node at (3.2,2.6) {\bf$ x_{3} ^{1}$};
				\node at (4.8,2.6) {\bf$ x_{3} ^{2}$};
				\node at (4.25,1.3) {\bf$ x_{2} ^{1}$};
				
				\node at (4,0.8) {$ K_1$};
				
				\node at (4,2.2) {\textbf{$ K_2$}};
				
				\node at (4,3.6) {\textbf{$ K_3$}};

				\node at (2.1,0.7) {\bf$ x_{12} ^{1}$};
				\node at (2.7,1.5) {\bf$ x_{12} ^{2}$};
				\node at (-1.2,4.3) {\bf$ x_{13} ^{1}$};
				\node at (1.2,4.3) {\bf$ x_{13} ^{3}$};
				\node at (-2,0.6) {\bf$ x_{23} ^{1}$};
				\node at (0,5.3) {\bf$ x_{13} ^{2}$};

				\draw (-2,1)--(1.9,0.9)--(-1,4)--(-2,1)--(0,5)--(-2,1)--(1,4)--(-2,1)--(2.4,1.6)--(-1,4);
				\draw (0,5)--(2.4,1.6)--(1,4)--(1.9,0.9)--(0,5)--(1,4);
				\draw (1,4)--(-1,4)--(0,5);
				\draw (1.9,0.9)--(2.4,1.6);

				\node at (0,0.2) { $H $};
				
				\node at (1.6,-0.6) { \tiny $ G(L^B)^{**}=H+K_1+K_2+K_3$};

			\end{scope}

		\end{tikzpicture}
		\caption{$\text{Illustration of Lemma \ref{blowl2}.}$}\label{figure4}
	\end{center}
\end{figure}

\begin{lemma}\label{blowl1}
	Suppose $L^B$ be a blow-up of a Boolean lattice $L\cong \mathbf{2}^n$ with $|[q_i]|\geq 2$ for atoms $q_i$, $1\leq i \leq n$ and $n\geq 3$. Then $G(L^B)^{**}=G(L^B)_{SR}$.
\end{lemma}

\begin{proof}
	First we show that  $V(G(L^B)^{**})=V(G(L^B))_{SR}$. 
	
	Let $x\in V(G(L^B))_{SR}$. Then there exists some $y\in V(G(L^B))_{SR}$ such that $x$ and $y$ are mutually maximally distant. Then by Lemma \ref{mutually}, gives that $x\wedge y \neq 0$, $x\nleq y$ and $y\nleq x$. Thus $x\in V(G(L^B)^{**})$. 
	
	Let $x\in V(G(L^B))^{**}=V(G(L^B))$.
	
	Suppose  $x\in Atoms(L^B)$. Then for some $y\in [x]$, $N(x)=N(y)$. This implies that $x$ and $y$ are mutually maximally distant. Hence, $x\in V(G(L^B)_{SR})$. 
	
	Now, suppose $x\notin Atoms(L^B)$. Then by Lemma \ref{distance}, there exist $y\in [x^*)$ such that $d(x,y)=3=diam(G(L^B))$. This shows that $x$ and $y$ are mutually maximally distant.


	Therefore $x\in V(G(L^B)_{SR})$. Hence $V(G(L^B)_{SR})=V(G(L^B))^{**}$.
	
	\par  Let $x$ be adjacent to $y$ in $G(L^B)^{**}$, that is,  either $([x]=[y])$ or $(a\wedge b \neq 0$, $a\nleq b$ and  $b\nleq a$ for every $a\in [x]$ and for every $b\in [y])$. Hence $[x] || [y]$, by Definition \ref{GLB}. 
	
	We show that $x$ is adjacent to $y$ in $G(L^B)_{SR}$.
	
	If $[x]=[y]$, then $N(x)=N(y)$.  Hence, $x$ and $y$ are mutually maximally distant. Therefore $x$ is adjacent to $y$ in $G(L^B)_{SR}$.
	
	Now, if $[x]\neq [y]$, then by adjacency of $x$ and $y$, we have $a\wedge b\neq 0$, $a\nleq b$ and $b\nleq a$, for every $a\in [x]$, $b\in [y]$.
	
	As $a \wedge b \neq 0$, $d(a,b)_{G(L^B)}\neq 1$. Therefore $d(a,b)_{G(L^B)}=\{2,3\}$. If $d(a,b)_{G(L^B)}=3=diam(G(L^B))$. Then  $a$ and $b$ are  mutually maximally distant for every $a\in [x]$ and $b\in [y]$. 
	Thus $x$ is adjacent to $y$ in $G(L^B)_{SR}$.		
	
	Now, suppose that $d(a,b)_{G(L^B)}=2
	$ and $c\in N_{G(L^B)}(a)$. Since $a\wedge c=0$ and $a\nleq b$, we claim that  $b\vee c\neq 1$. On the contrary assume that $b\vee c=1$, that is $1^*=(b \vee c)^*=b^*\wedge c^*=0$. This gives that $c^*\leq b^{**}$. Also,  $a\wedge c=0$ implies that $a\leq c^*$. This together gives that $a\leq b^{**}$. This means that $b^*\leq a^*$ and hence $[a]\leq [b]$, in particular, $[x]\leq [y]$, a contradiction to the fact that $[x] || [y]$.		 
	This means that $b\vee c\neq 1$.
	
	Now, if either  $b\wedge c=0$ or $b\wedge c\not =0$,  we have by Lemma \ref{alizadeh}, $d_{G(L^B)}(b,c)\leq 2$.
	Thus,  $d_{G(L^B)}(b,c)\leq d_{G(L^B)}(a,b)$.
	
	Similarly, we can show that $d_{G(L^B)}(a,e)\leq d_{G(L^B)}(a,b)=2$, for every $e\in N(b)$. This shows that $a$ and $b$ are mutually maximally distant. As $a\in[x]$ and $b\in [y]$, thus $x$ and $y$ are mutually maximally distant. Therefore $x$ is adjacent to $y$ in $G(L^B))_{SR}$. 
	
	 Let $x$ be adjacent to $y$ in $G(L^B))_{SR}$, that is, $x$ is mutually maximally distant with $y$ in $G(L^B)$. 	 By Lemma \ref{mutually},  $x\wedge y \neq 0$, $x || y$. Thus, $x$ and $y$ are adjacent in $G(L^B)^{**}$.
\end{proof}

\begin{note}\label{note3.22}
	If $|[q_i]|=1$ for some $i$, $1\leq i\leq n$, then by Lemma \ref{blowl2},  $K_{|[q_i]|}=K_1$. Further, by Lemma \ref{blowl1}, $G(L^B)^{*}=G(L^B)_{SR}$. Also, if $\#(|[q_i]|=1)=m$, then  $|V(G(L^B)_{SR})|=|Z^*(L^B)|-m$.
\end{note}

If $\Gamma(R)$ is the zero-divisor graph of a commutative ring $R$ with identity, then  the graph $\Gamma(R)^{* *}$ is defined  as follows: $V\left(\Gamma(R)^{* *}\right)=V(\Gamma(R))$ and two distinct vertices $a, b$ are adjacent in $\Gamma(R)^{* *}$ if and only if either $\operatorname{ann}(a)=\operatorname{ann}(b)$ or $a b \neq 0$ and $\operatorname{ann}(a b) \neq \operatorname{ann}(a) \cup \operatorname{ann}(b)$, where $\operatorname{ann}(a)=\{x \in R\;|\; xa=0\}$ (see \cite{nik}). Also, let $\Gamma(R)^*=\Gamma(R)$, if $\Gamma(R)$ is complete, otherwise, $\Gamma(R)^*$ is extracted from $\Gamma(R)^{* *}$ after deleting all isolated vertices.

A. Badawi \cite{badawi} introduced the annihilator graph of a commutative ring $R$ with identity denoted by $AG(R)$  with the vertex set $Z(R)^*=$ $Z(R) \backslash\{0\}$, and two distinct vertices $x$ and $y$ are adjacent if and only if $\operatorname{ann}(xy) \neq \operatorname{ann}(x) \cup \operatorname{ann}(y)$. It follows that each edge (path) of $\Gamma(R)$ is an edge (path) of $AG(R)$.

It is well known that every Boolean algebra or Boolean lattice $L\cong \mathbf{2}^n$ gives rise to a Boolean ring $R_L\cong \prod_1^n \mathbb{Z}_2$, and vice versa, with ring multiplication corresponding to the meet operation. Hence, we have $\operatorname{ann}(a)=a^\perp=(a^*]$, where $a^*$ is the pseudocomplement of $a$ in the Boolean lattice $L$.  Further, note that every element of a Boolean lattice $L$ has the unique complement which is also the pseudocomplement. Hence $a^{**}=a$ for every $a \in L$.

The following result is immediate from the above discussion.

\begin{lemma}\label{bool}
	The zero-divisor graph $G(L)$ of a Boolean lattice $L\cong \mathbf{2}^n$ is same as the zero-divisor graph $\Gamma(R_L)$ of a Boolean ring $R_L\cong \prod_1^n \mathbb{Z}_2$ derived from $L$. Hence $G(L)^{**}=\Gamma(R_L)^{**}$ and $G(L)_{SR}$ and  $\Gamma(R_L)_{SR}$.
\end{lemma}  

\begin{definition}[{Jejurkar and Joshi \cite[Definition 1.2]{rjj}}]\label{com}
	Let $L$ be a bounded lattice. The {\it comparability graph } of $L$ is an undirected, simple graph denoted by $Com(L)$, where the vertex set is $L\setminus \{0_L, 1_L\}$ and two vertices $a$ and $b$ are adjacent if and only if $a$ and $b$ are comparable. The complement of $Com(L)$  is the incomparability graph  $Incomp(L)$.
\end{definition}

\begin{lemma}\label{boolAGR}
	Let $L\cong \mathbf{2}^n$ be a Boolean lattice and $R_L\cong \prod_1^n \mathbb{Z}_2$ be a Boolean ring derived from $L$. Then $Incomp(L) = AG(R_L)$.
\end{lemma}	

\begin{proof}
	Clearly, $V(Incomp(L)) = V(AG(R_L))$. Let $a$ and $b$ two distinct adjacent vertices of $AG(R_L)$. Then  $\operatorname{ann}(a b) \neq \operatorname{ann}(a) \cup \operatorname{ann}(b)$. By \cite[Lemma 2.2]{nik0}, we have $\operatorname{ann}(a)\not\subseteq  \operatorname{ann}(b)$ and $\operatorname{ann}(b)\not\subseteq  \operatorname{ann}(a)$. This yields that $a^* \not\leq b^*$ and $b^* \not\leq a^*$ in $L$. Hence $a \not\leq b$ and $b \not\leq a$. Thus $a$ and $b$ are adjacent in $Incomp(G(L))$.
	
	Conversely, assume that $a$ and $b$ are adjacent in $Incomp(G(L))$. Then $a^* \not\leq b^*$ and $b^* \not\leq a^*$ which further gives $\operatorname{ann}(a)\not\subseteq  \operatorname{ann}(b)$ and $\operatorname{ann}(b)\not\subseteq  \operatorname{ann}(a)$. Again by \cite[Lemma 2.2]{nik0}, we have  $\operatorname{ann}(a b) \neq \operatorname{ann}(a) \cup \operatorname{ann}(b)$. 
\end{proof}

Let $L\cong \mathbf{2}^n$ be a Boolean lattice and $R_L\cong \prod_1^n \mathbb{Z}_2$ be a Boolean ring derived from $L$. 	Hence by Lemma \ref{boolAGR}, we have $E(\Gamma(R_L)^{**})= E(Incomp(G(L))) \cap E(\Gamma^c(R))$, where $\Gamma^c(R)$ is the complement of the zero-divisor graph $\Gamma(R)$. Since any two atoms of $L$ are not adjacent in $G^c(L)$ and consequently, in $\Gamma^c(R)$, it is clear that $V(\Gamma(R_L)^{**})$ will not contain atoms of $L$.

\begin{lemma}\label{mainlemma}
	
	Suppose $L^B$ is the blow-up of a Boolean lattice $L\cong \mathbf{2}^n$ with $n\geq 3$. Then the following statements hold.
	\begin{enumerate}
		\item If $L^B\cong L\cong \mathbf{2}^n$, then $\beta(G(L^B)_{SR})=n-2$.
		\item If $|[q_{i}]|\geq 2$ for every $i$,  $1\leq i\leq n$, then $\beta(G(L^B)_{SR})=2n-2$.
		\item If $\#(|[q_{i}]|=1)= m $  for some $i$, $1\leq i\leq n$, then $\beta(G(L^B)_{SR})=2n-m-2$.
	\end{enumerate}
\end{lemma}\label{blowl3}

\begin{proof}
	
		(1)  Follow from Lemma \ref{bool} and \cite[Lemma 3.2]{nik}.
		
		(2) By Lemma \ref{blowl2}, we have $G(L^B)_{SR}= H+ K_{|[q_{1}]|}+K_{|[q_{2}]|}+ \dots +K_{|[q_{n}]|}$. Also, it is well known that $\beta(K_{|[q_{1}]|}+K_{|[q_{2}]|}+ \dots +K_{|[q_{n}]|})=n$.  Since $G(L^B)^{**}=G(L^B)_{SR}$ and $G(L^B)^{**}$ has no isolated vertex. This implies that $\beta (G(L^B)_{SR})=\beta(H)+n$. It is enough to show that $\beta(H)=n-2$. Define a set  $A=\{[x]\mid x\in V(H)\}$. In $G(L^B)_{SR}$, $H([A])$ is a complete graph. Thus $\beta(H)=\beta(G(L^B)_{SR}[A])$. Note that from $(1)$, $G(L^B)_{SR}\cong H$ and hence  $\beta(H)=n-2$. Thus $\beta (G(L^B)_{SR})=\beta(H)+n=n-2+n=2n-2$.
		
		(3) Follow from $(1)$ and $(2)$. 
\end{proof}

Now, we are ready to state the main result of this paper.
\begin{theorem}\label{mainthm}
	Let $L^B$ be a blow-up of a Boolean lattice $L\cong \mathbf{2}^n$ $(n\geq 3)$ and $\#(|[q_i]|=1)=m$. Then $sdim_{M}(G(L^B))=|Z^*(L^B)|-2n+2$.	
\end{theorem}
\begin{proof}
	By Theorem \ref{galai} and Theorem \ref{gala}, $sdim_{M}(G(L^B))=\alpha(G(L^B))_{SR}$. Then by Lemma \ref{mainlemma} (3), we have $\beta(G(L^B))_{SR}=2n-m-2$, where $\#(|[q_i]|=1)=m$. Therefore, $sdim_{M}(G(L^B))=|V(G(L^B))_{SR}|-\beta(G(L^B))_{SR}= |V(G(L^B))_{SR}|-2n+m+2$. From Note \ref{note3.22}, we have    $sdim_{M}(G(L^B))=(|Z^*(L^B)|-m)-(2n-m-2)= |Z^*(L^B)|-2n+2$. Therefore, $sdim_{M}(G(L^B))=|Z^*(L^B)|-2n+2$.
\end{proof}
The following corollary immediately follows from  Theorem \ref{mainthm}. 

\begin{corollary}
	Let $L^B$ be a blow-up of a Boolean lattice $L\cong \mathbf{2}^n$ with $(n\geq 3)$. If $L^B\cong L\cong \mathbf{2}^n $, then $sdim_{M}(G(L^B))=2^n-2n$.
\end{corollary}


\section{Applications to graphs from Algebraic Structures}\label{applications}

In this section, we provide some applications of our results to the comaximal graph, the comaximal ideal graph, the zero-divisor graph of a reduced ring and the nonzero component graph of vector spaces.
\subsection{Comaximal graph of a ring}
\indent In \cite{pdsb}, Sharma and Bhatwadekar  introduced a graph $\Gamma_0(R)$ on a commutative ring $R$ with identity, whose vertices are the elements of $R$ and two distinct vertices $x$ and $y$ are adjacent if and only if $Rx + Ry = R$.  Maimani et al. \cite{hm} named the graph $\Gamma_0(R)$ studied by Sharma and Bhatwadekar as the {\it comaximal graph} of $R$. 

Maimani et al. \cite{hm} studied the subgraphs $\Gamma_1(R)$, $\Gamma_2(R)$ and $\Gamma_2^{\prime}(R)=\Gamma_2(R) \backslash J(R)$, where $\Gamma_1(R)$ is the subgraph of $\Gamma_0(R)$ induced on the set of units of $R$, $\Gamma_2(R)$ is the subgraph of $\Gamma_0(R)$ induced on the set of non-units of $R$ and $\Gamma_2^{\prime}(R)$ is the subgraph of $\Gamma_0(R)$ induced on the set of non-units of $R$ which are not in $J(R)$, the Jacobson radical of $R$, i.e., $\Gamma_2^{\prime}(R)=\Gamma_0(R) \backslash(U(R) \cup J(R))$. Moconja and Petrovi\'{c}
\cite{smm} shows that the comaximal graphs are blow-ups of Boolean graphs, the zero-divisor graphs of Boolean rings, equivalently, Boolean lattices. However, the construction of a Boolean lattice was not given.  The following result is essentially proved in \cite{pgg}.

\begin{theorem}[{Gadge et al. \cite[Theorem 3.16]{pgg}}]\label{comax}
	Let $R$ be a finite commutative ring with identity such that $|Max(R)|=n$. Then $\Gamma_2'(R)= G({L}^B)$, where $L^B$ is the blow-up of a Boolean lattice $L\cong \mathbf{2}^n$. 
	
\end{theorem}

The following result follows from  Theorem \ref{comax}, Theorem \ref{0-disblowup}, and Theorem \ref{mainthm}. 
\begin{theorem}
	Let $\Gamma_{2}'(R)$ be the comaximal graph of a commutative ring $R$ with identity and $|\mathrm{Max}(R)|=n$, $n\geq 3$. Then $sdim_{M}(\Gamma_{2}'(R))=|V(\Gamma_{2}'(R))|-2n+2$.
\end{theorem}

\subsection{Zero-divisor graph of a reduced ring}

Now, we compute the strong metric dimension of the zero-divisor graph of a reduced ring.

\begin{theorem}[{\cite[Remark 3.4]{jlkr}}, {\cite[Lemma 3.3]{sdj1}}]\label{zero}
	Let $\Gamma(R)$ be the ring-theoretic zero-divisor graph of a finite reduced commutative ring  $R$ with identity. Then $\Gamma(R)$ equals to the lattice-theoretic zero-divisor graph of $G(\prod_{i=1}^{n} C_i) $, where $C_i$'s are the chains with $|C_i|=|F_i|$, where $R=\prod_{i=1}^{n}F_i$ ($F_i$'s are finite fields.).
\end{theorem}

The following result follows from Theorem \ref{0-disblowup}, Theorem \ref{mainthm} and Theorem \ref{zero}.
\begin{corollary}[{R. Nikandish et al. \cite[Theorem 3.1]{nik}}]\label{reduced}
	Let $R$ be a ring. Then, the following hold.
	\begin{enumerate}
		
		\item If $R \cong \prod_1^n \mathrm{~F}_i$, where $F_i \neq \mathbb{Z}_2$ is a field for every $1 \leq i \leq n$, then $\operatorname{sdim}_M(\Gamma(R))=$ $\left|Z(R)^*\right|-2 n+2$.
		\item  If $R \cong \prod_1^n F_i \times \prod_1^m \mathbb{Z}_2$, where $F_i \neq \mathbb{Z}_2$ is a field for every $1 \leq i \leq n$, then $\operatorname{sdim}_ M(\Gamma(R))=\left|Z(R)^*\right|-2 n-2 m+2$.
	\end{enumerate}
\end{corollary}

\subsection{Comaximal ideal graph of a ring and co-annihilating ideal graph of  a ring}\label{coco}

Let $R$ be a commutative ring with identity and $Id(R)$ be the set of all ideals of $R$. Clearly, $Id(R)$ is a poset under set inclusion as a partial order. Then $(Id(R),\leq)$  is a modular, $1$-distributive lattice 
under the set inclusion as a partial order. Clearly, sup$\{I,J\}=I+J$ and  inf$\{ I, J\} = I\cap J$.
It is well known that the lattice $Id( R)$ is a complete lattice with the ideals $(0)$ and $R$ as its least and the greatest element, respectively.
Now, we denoted the lattice $Id(R)$ by $L$. Let $L^\partial$ be the dual of the lattice of $L$. Therefore in $L^\partial$, sup$_{L^\partial}\{I, J\} = I\cap J$ and
inf$_{L^\partial}\{I, J\} = I+ J$. The ideal $R$ is the least element of $L^\partial$, and the ideal $(0)$ is the greatest element of $L^\partial$. Further, by the duality, $   L^\partial$ is a $0$-distributive lattice. Moreover, the maximal ideals of $R$ are nothing but the atoms of $L^\partial$. Therefore, $L^\partial$ is an atomic lattice.

\par 	

\begin{definition}[{Ye and Wu \cite{mytw}, Akbari et al. \cite {saaa}}] Let $R$ be a commutative ring with identity. The {\it comaximal ideal graph}, $\mathbb{CG}(R)$ is a simple graph with its vertices the nonzero proper ideals of $R$ not contained in Jacobson radical $J(R)$ of $R$ and two distinct vertices $I$ and $J$ are adjacent if and only if $I+J=R$.

	The {\it co-annihilating ideal graph} of $R$, denoted by $\mathbb{CAG}(R)$ is a graph whose vertex set is the set of all nonzero proper ideals of $R$ and two distinct vertices $I$ and $J$ are adjacent whenever $ 
	\operatorname{ann}(I) \cap  \operatorname{ann} (J)=\{0\}$, where $\operatorname{ann}(I)=\{x\in R ~|~ xi=0~ \text{for all } i \in I\}$.
\end{definition}

In \cite{mytwl}, M. Ye et al. proved that the comaximal ideal graph $\mathbb{CG}(R)$ is the blow-up of the zero-divisor graph of a Boolean lattice $\mathbf{2}^n$. In fact, they proved,

\begin{theorem}[{M. Ye et al. \cite[Theorem 3.1]{mytwl}}]\label{comaxiid}
	Let $R$ be a ring with $|\text{Max}(R)|=n$, where $2 \leq n < \infty$. Then $\mathbb{CG}(R)$ is a blow-up of the zero-divisor graph of a Boolean lattice $\mathbf{2}^n$.
\end{theorem}

\begin{theorem}[{Khandekar and Joshi \cite[Theorem 5.1]{nkvj23}}]\label{comaxi}
	Let $R$ be a commutative ring with identity and let $Id(R)^\partial$ be the dual of the lattice $Id (R)$ of all ideals of $R$. Then $\mathbb{CG}(R)=G(Id(R)^\partial)$.
\end{theorem}

\begin{corollary}[{ \cite[Corollary 1.2]{saaa}}]\label{artinian}
	Let $R$ be an Artinian ring. Then $\mathbb{CAG}=\mathbb{CG}(R)$.
\end{corollary}
By  Theorem \ref{0-disblowup}, Theorem \ref{mainthm}, Theorem \ref{comaxiid},  and Theorem \ref{comaxi}, we have:
\begin{corollary} [{R. Shahriyari et al. \cite[Theorem 2.9]{smco}}]
	Let $R$ be a reduced commutative ring with identity, and $sdim_{SM}(\mathbb{CG}(R))$ is finite. Then, the following statements hold.  \begin{enumerate}
		\item If $|\text{Max}(R)|=2$, then $sdim_{M}(\mathbb{CG}(R))=sdim_{M}(\mathbb{CAG}(R))=1$.
		\item If $|\text{Max}(R)|=n\geq 3$, then $sdim_{M}(\mathbb{CG}(R))=sdim_{M}(\mathbb{CAG}(R))=2^n-2n$.
	\end{enumerate}
\end{corollary}

\begin{corollary} [{R. Shahriyari et al. \cite[Theorem 3.4]{smco}}]
	Suppose that $R\cong R_1 \times R_2 \times \dots \times R_n$, where $R_i$ is an Artinian local ring and $|A(R_i)^{*}|\geq 1$, for every $1\leq i\leq n$. Then $sdim_{M}(\mathbb{CG}(R))=sdim_{M}(\mathbb{CAG}(R))=|V(\mathbb{CG}(R)))|-2n+2$.
\end{corollary}

\subsection{Component graphs of vector spaces}
\par Angsuman Das \cite{das2} defined and studied the nonzero component graph union graph of a finite-dimensional vector space. Let $\mathbb{V}$ be a vector space over field $\mathbb{F}$ with  $\mathcal{B}=\{v_1,\dots,v_n\}$ as a basis and $0$ as the null vector. Then, any vector $a\in \mathbb{V}$ can be uniquely expressed in the linear combination of the form $a=a_1v_1+\dots+a_nv_n$. We denote this representation as a basic representation of $a$ with respect to $\{v_1,\dots,v_n\}$. Define the skeleton of $a$ with respect to $\mathcal{B}$, as
\[S_{\mathcal{B}}(a)=\{v_i\mid a_i\neq 0,  a=a_1v_1+\dots+a_nv_n\}.\] 

Angsuman Das \cite{das2}  defined the \textit{nonzero component union  graph} $\mathbb{UG(V)}$ with respect to $\mathcal{B}$ as follows: The vertex set of graph  $\mathbb{UG(V)}$ is $\mathbb{V}\setminus \{0\}$ and for any $a,b\in \mathbb{V}\setminus \{0\}$, $a$ is adjacent to  $b$ if and only if $S_{\mathcal{B}}(a)\cup S_{\mathcal{B}}(b)=\mathcal{B}$.

\par In \cite{ncv}, Khandekar et al. gave a relation between the skeleton union graph of a finite-dimensional vector space and the zero-divisor graph of the blow-up of a Boolean lattice. Hence, we have the following result.

\begin{theorem}[Khandekar et al.\cite{ncv}]
	Let $\mathbb{V}$ be a $n$-dimensional vector space over a field $\mathbb{F}$. Then $\mathbb{UG(V)}=G({L}^B)$ $\vee K_t$, where $t=|V_{12\dots n}|=(|\mathbb{F}|-1)^n$ and $L^B$ is the blow-up of a Boolean lattice $L\cong \mathbf{2}^n$.
\end{theorem}
By Theorem \ref{mainthm}, we have the following result.

\begin{theorem}
	Let $\mathbb{UG(V)}$ be the nonzero component union graph of vector spaces with $dim(\mathbb{V})=n\geq 3$. Then $sdim_{M}(\mathbb{UG(V)})=|V(\mathbb{UG(V)})|-n+2$.
\end{theorem}

\par\noindent 

\noindent\textbf{Funding:}\\
First author: None.\\
Second author: Supported by DST(SERB) under the scheme CRG/2022/002184.

\noindent\textbf{Conflict of interest:} The authors declare that there is no conflict of
interests regarding the publishing of this paper.

\noindent\textbf{Authorship Contributions :} Both authors contributed to the study on the strong metric dimension of a zero-divisor graph of a poset. Both authors read and approved the final version of the manuscript.

\noindent \textbf{Data Availability Statement :} Data sharing does not apply to this article, as no datasets were generated or analyzed during the current study.

\end{document}